# Vector invariants in arbitrary characteristic

by

Frank D. Grosshans
Department of Mathematics,
West Chester University of Pennsylvania
West Chester, PA 19383
fgrosshans@wcupa.edu

**Abstract.** Let $k$ be an algebraically closed field of characteristic $p \geq 0$. Let $H$ be a subgroup of $GL_n(k)$. We are interested in the determination of the vector invariants of $H$. When the characteristic of $k$ is 0, it is known that the invariants of $d$ vectors, $d \geq n$, are obtained from those of $n$ vectors by polarization. This result is not true when char $k = p > 0$ even in the case where $H$ is a torus. However, we show that the algebra of invariants is always integral over the algebra of polarized invariants and when $H$ is reductive is actually the $p$ - root closure of that algebra. We also give conditions for the algebras to be equal, relating equality to good filtrations and saturated subgroups. We conclude with examples where $H$ is finite or a classical group or is a certain kind of unipotent subgroup of $GL_n(k)$.

## §1. Introduction

Let $k$ be an algebraically closed field of characteristic $p \geq 0$. Let $H$ be a subgroup of $GL_n = GL_n(k)$. We shall be interested in the determination of the vector invariants of $H$. When the characteristic of $k$ is 0, this determination is greatly simplified by the following Theorem [W, p.44]: all invariants of $d$ vectors, $d \geq n$, are known once one knows the invariants for $n$ vectors; indeed, the invariants of $d$ vectors are obtained from those of $n$ vectors by polarization. This theorem is not true when char $k = p > 0$ even in the case where $H$ is a torus. The purpose of this paper is to find the correct statement of the Theorem when char $k = p > 0$. Before doing that, we need to introduce some notation which will be used throughout this paper.

**Notation.**

$k$: algebraically closed field of characteristic $p \geq 0$.

$M_{n,d}$: the algebra of $n \times d$ matrices over $k$ where $d \geq n$. If $x \in M_{n,d}$, we shall write $x = (x_1, \ldots, x_d)$
    where $x_i$ is the $i$th column of $x$.

$\Delta(x_1, \ldots, x_d)$ = determinant of $(x_1, \ldots, x_n)$.

$H$ : subgroup of $GL_n$.

$Z$: affine variety on which $H$ acts regularly via a mapping $H \times Z \to Z$ denoted by $(h, z) \to hz$. (In what follows, the variety $Z$ could be $\{0\}$.)

$j: k[M_{n,n} \times Z] \to k[M_{n,d} \times Z]$ given by $(jF)(x_1, \ldots, x_d; z) = F(x_1, \ldots, x_n; z)$ for $F \in k[M_{n,n} \times Z], z \in Z$ and $x = (x_1, \ldots, x_d) \in M_{n,d}$.

The subgroup $H$ acts on $M_{n,d} \times Z$ by $h(x; z) = (hx; hz)$ where $hx$ is the ordinary matrix product of $h$ and $x$. This action gives an action of $H$ on $k[M_{n,d} \times Z]$ by $(hF)(x; z) = F(h^{-1}x; h^{-1}z)$ for $F \in k[M_{n,d} \times Z]$. Let $^H k[M_{n,d} \times Z] = \{F \in k[M_{n,d} \times Z] : hF = F \text{ for all } h \in H\}$. The mapping $j$ is $H$-equivariant so if $F \in {}^H k[M_{n,n} \times Z]$, then $jF \in {}^H k[M_{n,d} \times Z]$.

The group $GL_d$ acts on $M_{n,d} \times Z$ by $g * (x; z) = (xg^{-1}; z)$ where $xg^{-1}$ is the ordinary matrix product of $x$ and $g^{-1}$. This gives an action of $GL_d$ on $k[M_{n,d} \times Z]$ by $(g * F)(x; z) = F(xg; z)$ for $F \in k[M_{n,d} \times Z]$. If $F \in k[M_{n,n} \times Z]$, then $(g * jF)(x; z) = (jF)(xg; z)$. The actions of $H$ and $GL_d$ commute so $GL_d$ sends $^H k[M_{n,d} \times Z]$ to itself. Let $GL_d * j(^H k[M_{n,n} \times Z])$ be the algebra generated by all $g * jf$ for $g \in GL_d, f \in {}^H k[M_{n,n} \times Z]$. This is the algebra of invariants of $d$ vectors obtained from those of $n$ vectors by polarization.

Three subgroups of $GL_d$ will occur frequently: the subgroup $B$ which consists of upper triangular matrices, the subgroup $U$ which consists of those elements in $B$ with 1's on the diagonal and the subgroup $T$ which consists of all diagonal matrices in $B$.

**Definition.** Let char $k = p > 0$ and let $R$ and $S$ be $k$-algebras with $R \subset S$. We say that $S$ is *contained in the $p$-root closure of $R$* (or *purely inseparable closure of $R$*) if for every $s \in S$, there is a non-negative integer $m$ so that $s^{p^m} \in R$. (If char $k = 0$, $S$ is contained in the $p$-root closure of $R$ if and only if $R = S$.)

Given this terminology, we ask: what is the relationship between $^H k[M_{n,d}]$ and $GL_d * j(^H k[M_{n,n}])$? In §2, we show that $^H k[M_{n,d}]$ is integral over $GL_d * j(^H k[M_{n,n}])$. In §3, we refine this result and show that $^Hk[M_{n,d}]$ is contained in the $p$-root closure of $GL_d * j(^Hk[M_{n,n}])$ when $H$ is a "Grosshans subgroup" of $GL_n$. (The class of Grosshans subgroups includes all finite subgroups, all reductive subgroups and their maximal unipotent subgroups.) Next, we consider the question as to when $^Hk[M_{n,d}] = GL_d * j(^Hk[M_{n,n}])$. In §4, we show that equality holds if $GL_d * j(^Hk[M_{n,n}])$ has a good $GL_d$ filtration. In §5, we obtain a partial converse: if $H$ is a saturated subgroup of $GL_n$ and if $GL_d * j(^Hk[M_{n,n}]) = {}^Hk[M_{n,d}]$, then $GL_d * j(^Hk[M_{n,n}])$ has a good $GL_d$ filtration. In §6, we study examples where $H$ is finite or a classical group or is a certain kind of unipotent subgroup of $GL_n$.

## §2. Integral closure

For a moment, let $G$ be any reductive group. Let $B = TU$ be a Borel subgroup of $G$ where $U$ is a maximal unipotent subgroup and $T$ is a maximal torus normalizing $U$. Let $A$ be a commutative $k$-algebra on which $G$ acts rationally via a mapping $G \times A \to A$ denoted by $(g, a) \to g \cdot a$. Let $G \cdot A^U$ be the algebra generated by all $g \cdot f$ for $g \in G, f \in A^U$.

**Theorem 1.** [G2; Theorem 5, p.128] The algebra $A$ is integral over $G \cdot A^U$.

We apply this theorem to the case where the group $G$ is $GL_d$.

**Lemma 2.** Let $f \in k[M_{n,d}]^U$. Then $f \in j(k[M_{n,n}])$. In particular, $^Hk[M_{n,d}]^U = (GL_d * j(^Hk[M_{n,n}]))^U$.
**Proof.** Let $f' \in k[M_{n,n}]$ be defined by $f'(x_1, \ldots, x_n) = f(x_1, \ldots, x_n, 0, \ldots, 0)$. We shall show that $(jf')(x) = f(x)$ for all $x = (x_1, \ldots, x_d) \in M_{n,d}$. In doing this, we may assume that the set $\{x_1, \ldots, x_n\}$ is linearly independent over $k$ since the set of such $x$ is Zariski dense in $M_{n,d}$. Given such an $x$, we define an element $u \in U$ as follows. When $j \leq n$, the $j$th column of $u$ has 1 in position $(j, j)$ and 0's elsewhere. When $j > n$, we find scalars $u_{1j}, \ldots, u_{nj}$ in $k$ so that $x_j = -(u_{1j}x_1 + \ldots + u_{nj}x_n)$. Then, the first $n$ entries in the $j$th column of $u$ are $u_{1j}, \ldots, u_{nj}$ and all the other entries above the diagonal in the $j$th column are 0. The matrix $xu$ has $j$th column $x_j + u_{1j}x_1 + \ldots + u_{nj}x_n = 0$. Hence, $xu = (x_1, \ldots, x_n, 0, \ldots, 0)$. Since $f$ is $U$-invariant, we have $(jf')(x) = (jf')(x_1, \ldots, x_d)$

$$= f'(x_1, \ldots, x_n) = f(x_1, \ldots, x_n, 0, \ldots, 0) = f(xu) = (u * f)(x) = f(x).$$

This proves the first statement. Since the mapping $j$ is $H$-equivariant, we then have ${}^H k[M_{n,d}]^U \subset (j({}^H k[M_{n,n}]))^U \subset (GL_d * j({}^H k[M_{n,n}]))^U$. As for the reverse inclusion, since $(GL_d * j({}^H k[M_{n,n}])) \subset {}^H k[M_{n,d}]$, we have $(GL_d * j({}^H k[M_{n,n}]))^U \subset {}^H k[M_{n,d}]^U$. □

**Note 1.** The algebra $k[M_{n,d}]^U$ can be described explicitly, see for example [G1; Theorem 13.3, p. 78]. It is generated by minors whose rows can be chosen arbitrarily but whose columns must be chosen to have the form $1, 2, \ldots, r$.

**Theorem 3.** ${}^H k[M_{n,d} \times Z]$ is integral over $GL_d * j({}^H k[M_{n,n} \times Z])$. In particular, ${}^H k[M_{n,d}]$ is integral over $GL_d * j({}^H k[M_{n,n}])$.

**Proof.** We have $k[M_{n,d} \times Z] = k[M_{n,d}] \otimes k[Z]$. Now, $GL_d$ fixes $k[Z]$ (since $GL_d$ acts trivially on $Z$). Hence, by Lemma 2, we have $k[M_{n,d} \times Z]^U = \{f \in k[M_{n,d} \times Z] : u * f = f \text{ for all } u \in U\}$

$$= k[M_{n,d}]^U \otimes k[Z] \subset j(k[M_{n,n}]) \otimes k[Z] = j(k[M_{n,n} \times Z]).$$

Therefore, ${}^H k[M_{n,d} \times Z]^U \subset {}^H j(k[M_{n,n} \times Z]) = j({}^H k[M_{n,n} \times Z])$ since the mapping $j$ is an $H$-equivariant injection. We may now apply Theorem 1 with $G = GL_d$ and $A = {}^H k[M_{n,d} \times Z]$. □

**Note 2:** subgroups of $SL_n$. Let $j: k[M_{n,n-1}] \to k[M_{n,d}]$ given by $(jF)(x_1, \ldots, x_d) = F(x_1, \ldots, x_{n-1})$ for $F \in k[M_{n,n-1}]$. Let $H$ be a subgroup of $SL_n$. According to the Note 1, $k[M_{n,d}]^U \subset j(k[M_{n,n-1}])[\Delta]$. This allows us to extend Theorem 3 as follows: ${}^H k[M_{n,d} \times Z]$ is integral over $GL_d * j({}^H k[M_{n,n-1} \times Z])[\Delta]$. In particular, ${}^H k[M_{n,d}]$ is integral over $GL_d * j({}^H k[M_{n,n-1}])[\Delta]$.

**Theorem 4.** Let $H$ be a subgroup of $SL_n$ and $F \in {}^H k[M_{n,d}]$. There is a non-negative integer $e$ so that $\Delta^e F \in GL_d * j({}^H k[M_{n,n}])$. In particular, the quotient field of $GL_d * j({}^H k[M_{n,n}])$ is equal to the quotient field of ${}^H k[M_{n,d}]$.

**Proof.** Let $d = n + r$. The set $GL_n \times M_{n,r} = \{(x_1, \ldots, x_d) \in M_{n,d} : \Delta(x_1, \ldots, x_d) \neq 0\}$ is open in $M_{n,d}$. Hence, we may assume that $F \in {}^H k[GL_n \times M_{n,r}]$. We shall denote an element in $GL_n \times M_{n,r}$ by $(g, x)$ where $g \in GL_n$ and $x \in M_{n,r}$.

The algebras ${}^H k[GL_n \times M_{n,r}]$ and $k[M_{n,r} \times H\backslash GL_n]$ are isomorphic. Indeed, given $f \in$

$k[M_{n,r} \times H\backslash GL_n]$, we define $\varphi f \in {}^H k[GL_n \times M_{n,r}]$ by $(\varphi f)(g, x) = f(g^{-1}x, Hg)$. Also, given $f' \in {}^H k[GL_n \times M_{n,r}]$, we define $\psi f' \in k[M_{n,r} \times H\backslash GL_n]$ by $(\psi f')(x, Hg) = f'(g, gx)$. From this we see that for $F \in {}^H k[GL_n \times M_{n,r}]$, $(\psi F)(x, Hg) = \sum a_i(x) b_i(Hg)$ where $a_i \in k[M_{n,r}]$ and $b_i \in k[H\backslash GL_n] = {}^H k[GL_n] = {}^H k[M_{n,n}][1/\Delta]$. Then $F(g, x) = \varphi(\psi F)(g, x) = (\psi F)(g^{-1}x, Hg) = \sum a_i(g^{-1}x) b_i(Hg)$.

Let $v$ be any column of $x$ and let $g \in GL_n$. The vector $w = g^{-1}v$ satisfies the equation $gw = v$. Solving this system of equations by Cramer's rule, we see that the $i$th coordinate of $w$ is the quotient of two determinants, namely $\det(g_1, \ldots, g_{i-1}, v, g_{i+1}, \ldots, g_n)/\det(g_1, \ldots, g_n)$ where $g_j$ is the $j$th column of $g$. For $j = r+1, \ldots, d$ and $i = 1, \ldots, n$, let $\Delta_{ij}$ be the determinant of the matrix $(g_1, \ldots, g_{i-1}, v_j, g_{i+1}, \ldots, g_n)$. It follows that each $a_i(g^{-1}x)$ is a polynomial in $\Delta_{ij}/\Delta$. The polynomial $\Delta$ is in ${}^H k[M_{n,n}]$ since $H$ is a subgroup of $SL_n$. Also, each $\Delta_{ij}$ is in $GL_d * j({}^H k[M_{n,n}])$ since it has the form $g * \Delta$ where $g$ is a suitable permutation matrix in $GL_d$. □

**Corollary 5.** Let $H$ be a subgroup of $SL_n$. Then ${}^H k[M_{n,d}]$ is the integral closure of $GL_d * j({}^H k[M_{n,n}])$ in its quotient field.
**Proof.** Apply Theorems 3 and 4. □

**Example:** $H$ torus, ${}^H k[M_{n,d}] \supsetneq GL_d * j({}^H k[M_{n,n}])$. Let char $k = p > 0$, $n = 2$ and $d = 2p$. We shall denote the entries in the first row of $m \in M_{2,2p}$ by $x_1, \ldots, x_{2p}$ and in the second row by $y_1, \ldots, y_{2p}$. If $g \in GL_{2p}$, then $g * x_r = \sum x_s g_{sr}$. In particular, if $i \geq p$, then

(*) $\quad (g * x_r)^i = (g * x_r)^p (g * x_r)^{i-p} = (\sum x_s g_{sr})^p (\sum x_s g_{sr})^{i-p} = (\sum x_s^p g_{sr}^p)(\sum x_s g_{sr})^{i-p}$.

Now, let $H$ be the subgroup of $GL_2$ consisting of diagonal matrices with entries on the diagonal $a^{-1}$, $a^{2p}$ where $a \in k^*$. It may be shown that ${}^H k[M_{2,2}]$ is generated by monomials $x_1^i x_2^j y_h$ where $h = 1$ or $2$ and $i + j = 2p$. Since $i + j = 2p$, either $i \geq p$ or $j \geq p$. Therefore, in each monomial that appears in $g * x_1^i x_2^j y_h$, some $x_r$ has exponent at least $p$ by (*). It follows that the polynomial $F = x_1 x_2 \ldots x_{2p} y_{2p}$ is in ${}^H k[M_{2,2p}]$ but is not in $GL_d * j({}^H k[M_{n,n}])$.

## §3. $p$ - closure

**Theorem 6** [vdK1; Sublemma A.5.1, p.88 or vdK2]. Suppose that char $k = p > 0$. Let $X$ and $Y$

be affine varieties and let $f: X \to Y$ be a proper bijective morphism. Then $k[X]$ is contained in the $p$-root closure of $k[Y]$.

**Theorem 7.** Let $(x; z)$ and $(y; z') \in M_{n,d} \times Z$. If there is an $F \in {}^Hk[M_{n,d} \times Z]$ so that $F(x; z) \neq F(y; z')$, then there is an $f \in GL_d * j({}^Hk[M_{n,n} \times Z])$ so that $f(x; z) \neq f(y; z')$.

**Proof.** We proceed by induction on $d$. When $d = n$, the two algebras ${}^Hk[M_{n,d} \times Z]$ and $GL_d * j({}^Hk[M_{n,n} \times Z])$ are equal. So let us assume that $d > n$ and let $x = (x_1, \ldots, x_d)$. Since $d > n$, we may assume (perhaps after renumbering) that $x_d = -(a_1 x_1 + \ldots + a_{d-1} x_{d-1})$. Let $u \in U$ be the matrix all of whose entries above the diagonal are 0 except for the $d$th column which has entries $a_1, \ldots a_{d-1}, 1$. Then $xu = (x_1, \ldots, x_{d-1}, 0)$. Both ${}^Hk[M_{n,d} \times Z]$ and $GL_d * j({}^Hk[M_{n,n} \times Z])$ are invariant under $GL_d$. Hence, we may replace $x$ by $xu$, $y$ by $yu$, $F$ by $u^{-1} * F$ and assume from the outset that $(x; z) = (x_1, \ldots, x_{d-1}, 0; z)$. We may also assume that $F$ is a $T$-weight vector with weight $\chi = e_1 \chi_1 + \ldots + e_d \chi_d$ where $e_i$ is the ordinary polynomial degree of $F$ in the variables corresponding to the $i$th column of an $n \times d$ matrix and $\chi_i(\text{diag}(t_1, \ldots, t_d)) = t_i$.

If $e_d = 0$, then $F$ only depends on the first $d - 1$ columns and we may apply the induction hypothesis. Otherwise, $e_d > 0$ and, then, $F(x; z) = 0$ since $x_d = 0$. Now, according to Theorem 3, $F$ is integral over $GL_d * j({}^Hk[M_{n,n} \times Z])$. Hence, there is a positive integer $m$ and $r_i \in GL_d * j({}^Hk[M_{n,n} \times Z])$ for $i = 0, \ldots, m-1$ so that
(*)  $$F^m + r_{m-1} F^{m-1} + \ldots + r_0 = 0.$$
Since $F$ has $T$-weight $\chi$, we may assume that $r_i$ has $T$-weight $(m-i)\chi$. If some $r_i(x; z) \neq r_i(y; z')$, we are finished. Otherwise, $r_i(y; z') = r_i(x; z) = 0$ since $r_i$ has $T$-weight $(m-i)\chi$ and $e_d > 0$. Evaluating (*) at $(y; z')$ gives $F^m(y; z') = 0$ so $F(y; z') = 0 = F(x; z)$. This contradicts our choice of $F$. □

**Theorem 8.** If $H$ is a reductive (not necessarily connected) subgroup of $GL_n$, then ${}^Hk[M_{n,d} \times Z]$ is contained in the $p$-root closure of $GL_d * j({}^Hk[M_{n,n} \times Z])$.

**Proof.** Let $X$ (resp. $Y$) be the affine variety such that $k[X] = {}^Hk[M_{n,d} \times Z]$ (resp. $k[Y] = GL_d * j({}^Hk[M_{n,n} \times Z])$). Let $\pi: X \to Y$ be the map corresponding to the inclusion, $\pi^*$, of $k[Y]$ in $k[X]$. According to Theorem 3, $\pi$ is finite and, so is proper and surjective. Let $\mu: M_{n,d} \times Z \to X$ correspond to the inclusion, $\mu^*$, of $k[X]$ in $k[M_{n,d} \times Z]$. It is known that $\mu$ is surjective [N; Theorem 3.5, p.61]. Then, the mapping $\pi$ is injective. Indeed, let $x'$ and $y'$ be distinct points in $X$. Since $\mu$ is surjective, we can find points $x$ and $y$ in $M_{n,d} \times Z$ such that $\mu(x) = x'$ and $\mu(y) = y'$. Since $x'$ and $y'$ are distinct, there is a function $F \in k[X]$ so that $(\mu^*F)(x) \neq (\mu^*F)(y)$. By Theorem 7, there is a function $f \in k[Y]$ so that $(\mu^*\pi^*f)(x) \neq (\mu^*\pi^*f)(y)$, i.e., $f(\pi(x')) \neq f(\pi(y'))$. We may now apply Theorem 6 to the mapping $\pi: X \to Y$. $\square$

**Corollary 9.** If $H$ is a reductive (not necessarily connected) subgroup of $GL_n$, then $^Hk[M_{n,d}]$ is contained in the $p$ - root closure of $GL_d * j(^Hk[M_{n,n}])$. If $H$ is a reductive subgroup of $SL_n$, then $^Hk[M_{n,d}]$ is the $p$ - root closure of $GL_d * j(^Hk[M_{n,n}])$ in its quotient field.

**Proof.** For first statement, let $Z = \{0\}$ in Theorem 8. For the second, apply Theorem 4. $\square$

**Theorem 10.** If $H$ is a Grosshans subgroup of $GL_n$, then $^Hk[M_{n,d}]$ is contained in the $p$ - root closure of $GL_d * j(^Hk[M_{n,n}])$. If $H$ is a Grosshans subgroup of $SL_n$, then $^Hk[M_{n,d}]$ is the $p$ - root closure of $GL_d * j(^Hk[M_{n,n}])$ in its quotient field.

**Proof.** In this proof we denote $GL_n$ by $G$. There is an affine $G$ - variety $Z$ and a point $z_0 \in Z$ so that $k[Z] = k[G/H] = k[G]^H$ and $G/H$ is isomorphic to the orbit $Gz_0$; furthermore, the mapping $\rho$: $^Gk[M_{n,d} \times Z] \to {}^Hk[M_{n,d}]$ given by $(\rho F)(x) = F(x; z_0)$ is a surjective isomorphism [G1, Theorems 4.3, p.20 and 9.1, p. 49 or Exercise 3, p. 53].

If $F \in {}^Gk[M_{n,n} \times Z]$, then $\rho(g * jF) = g * \rho(jF)$ for all $g \in GL_d$. Indeed, for $x \in M_{n,d}$, we have

$$\rho(g * jF)(x) = (g * jF)(x; z_0) = (jF)(xg; z_0) = \rho(jF)(xg) = (g * \rho(jF))(x).$$

Also, $\rho(jF) = j(\rho F)$ since for $x = (x_1, \ldots, x_d) \in M_{n,d}$, we have $\rho(jF)(x_1, \ldots, x_d)$
$= (jF)(x_1, \ldots, x_d; z_0) = F(x_1, \ldots, x_n; z_0) = (\rho F)(x_1, \ldots, x_n) = (j(\rho F))(x_1, \ldots, x_d)$.

Now, let $f \in {}^Hk[M_{n,d}]$. Then, $f = \rho(F)$ for some $F \in {}^Gk[M_{n,d} \times Z]$. By Theorem 8 with $H = G$, we can find a non - negative integer $m$ so that $F^{p^m}$ is a sum of terms, each having the form $c(g_1 * jF_1) \ldots (g_r * jF_r)$ where $g_i \in GL_d$ and $F_i \in {}^Gk[M_{n,n} \times Z]$. By what we proved above, $\rho(g_i * jF_i) = g_i * \rho(jF_i) = g_i * j(\rho F_i) \in GL_d * j(^Hk[M_{n,n}])$. Since $f^{p^m} = \rho(F^{p^m})$, the first statement in the Theorem is proved. The second follows from Theorem 4. □

### §4. Good filtrations

We begin this section by recalling some concepts valid for any connected reductive group $G$ [G1; §15, p. 86 or G2; §2, p.129]. We also review some facts related to "good filtrations". Later on, we apply these theorems to the case where $G = \mathrm{GL}_d$ acts on $k[M_{n,d}]$.

**The homomorphism $h$.** Let $G$ be a connected reductive group. The action of $G$ on $k[G]$ by right (resp. left) translation will be denoted by $r_g$ (resp. $\ell_g$). Thus, for $g, x \in G$ and $f \in k[G]$, we have $(r_g f)(x) = f(xg)$ (resp. $(\ell_g f)(x) = f(g^{-1}x)$). Let $T$ be a maximal torus in $G$ and let $X(T)$ be its character group. Let ( , ) be a positive definite symmetric bilinear form on $X(T) \otimes_{\mathbf{Z}} \mathbb{R}$ which is invariant under the Weyl group of $T$. Let $B$ be a Borel subgroup of $G$ which contains $T$ and has maximal unipotent subgroup $U$. Let $\Phi^+$ be the set of roots of $B$ relative to $T$. A partial order may be defined on $X(T)$ by $\chi' > \chi$ if $\chi' - \chi$ is a sum of roots in $\Phi^+$. Let $X^+(T)$ be the set of all dominant weights. We define a homomorphism $h: X(T) \to \mathbf{Z}$ by $h(\chi) = \Sigma 2(\chi, \alpha)/(\alpha, \alpha)$ where the sum is over all $\alpha \in \Phi^+$. The mapping $h$ satisfies the following conditions: (i) $h(\omega)$ is a non-negative integer for every $\omega \in X^+(T)$; (ii) whenever $\chi, \chi' \in X^+(T)$ satisfy $\chi' > \chi$, then $h(\chi') > h(\chi)$. Let $B^- = TU^-$ be the Borel subgroup of $G$ which is "opposite" to $B$. For $\omega \in X^+(T)$, we define subspaces $Y(\omega)$ of $k[G/U^-]$ as follows: $Y(\omega) = \{f \in k[G/U^-] : r_t f = \omega(t)^{-1} f \text{ for all } t \in T\}$. The group $G$ acts on $Y(\omega)$ by left translation.

**Good filtrations.** Let $V$ be a rational $G$-module. A *good filtration* of $V$ is given by a sequence of $G$-stable subspaces $0 = V_0 \subset V_1 \subset \ldots$ such that $V = \cup V_i$ and each $V_i/V_{i-1}$ is isomorphic to some $Y(\omega_i)$, $\omega_i \in X^+(T)$. If char $k = 0$, then all such $Z$ have good filtrations since actions of $G$ are completely reducible. If char $k > 0$, then locally finite $G$-modules may or may not have good filtrations.

(D1) [D3; 12.1.6, p.236] A subset $\pi$ of $X^+(T)$ is called *saturated* if whenever $\mu \in X^+(T)$ and $\lambda$

$\in \pi$ with $\mu \leq \lambda$, then $\mu \in \pi$. Let $\pi$ be saturated subset of $X^+(T)$ and let $V$ be a $G$-module having a good filtration. Let $O_\pi(V)$ be the unique $G$-submodule of $V$ which is maximal among those belonging to $\pi$, i.e., each highest weight in $O_\pi(V)$ is in $\pi$. Then $O_\pi(V)$ and $V/O_\pi(V)$ have good filtrations.

(D2) [D3; Proposition 3.2.4, p.36] Suppose that $V$ is a $G$-module with a good filtration and that $V'$ is a $G$-submodule which also has a good filtration. Then $V/V'$ has a good filtration.

**The algebras gr$A$ and hull$_\nabla$(gr$A$).** Let $A$ be a commutative $k$-algebra on which $G$ acts rationally via a mapping $G \times A \to A$ denoted by $(g, a) \to g \cdot a$. For $a \in A$, the linear span of all the elements $g \cdot a$ will be denoted by $<G \cdot a>$. Since $G$ acts rationally, each such vector space is finite-dimensional. We next construct a graded algebra, gr$A$, associated to this action. For each positive integer $m$, let $A_m = \{a \in A : h(\chi) \leq m$ for all weights $\chi$ of $T$ on $<G \cdot a>\}$; we set $A_{-1} = \{0\}$. As a vector space, gr$A$ is the direct sum of all the spaces $A_m/A_{m-1}$, $m = 0,1,...$ . If $a \in A_q$ and $b \in A_r$, we define $(a + A_{q-1})(b + A_{r-1}) = ab + A_{q+r-1}$. In this way, gr$A$ becomes a commutative $k$-algebra on which $G$ acts rationally. There is an isomorphism between the algebras $(\text{gr}A)^U$ and $A^U$ which is $T$-equivariant and extends to a $T$-equivariant algebra homomorphism $\varphi:$ gr$A \to A^U$ such that $\varphi(u \cdot v) = \varphi(v)$ for all $u \in U^-$ and $v \in$ gr$A$. Indeed, let $v \in A_m/A_{m-1}$ be written as a sum of $T$-weight vectors: $(a + A_{m-1}) + \sum(b_i + A_{m-1})$ where $a$ has $T$-weight $\chi$ with $h(\chi) = m$ and $b_i$ has $T$-weight $\mu_i$ with $h(\mu_i) < m$. Then $a \in A^U$ and $\varphi(v) = a$.

Since $T$ normalizes $U$ and $U^-$, $A^U$ is invariant under the action of $T$ and $T$ acts on $k[G/U^-]$ by right translation. With respect to these actions, let hull$_\nabla$(gr$A$) = $(A^U \otimes_k k[G/U^-])^T$. The group $G$ acts on hull$_\nabla$(gr$A$) by the trivial action on $A^U$ and by left translation on $k[G/U^-]$. If we choose a basis $\{a_i\}$ of $A^U$ such that $t \cdot a_i = \omega_i(t)a_i$ for all $t \in T$, then we see that hull$_\nabla$(gr$A$) is a direct sum of the subspaces $a_i \otimes Y(\omega_i)$.

A morphism from $G/U^-$ to $A^U$ is a mapping $\alpha$ of the form $\alpha(gU^-) = \Sigma_i f_i(g)b_i$ where each $f_i \in k[G/U^-]$, $b_i \in A^U$ and where only finitely many of the functions $f_i$ are non-zero. The group $G$ acts by left translation, namely, $(g \cdot \alpha)(g_1 U^-) = \alpha(g^{-1}g_1 U^-) = \Sigma_i (\ell_g f_i)(g_1)b_i$. The vector space $A^U \otimes k[G/U^-]$ may be identified with the vector space consisting of all morphisms from $G/U^-$ to $A^U$ by associating to $\Sigma_i b_i \otimes f_i$ the morphism which sends $gU^-$ to $\Sigma_i f_i(gU^-)b_i$. The morphisms associated to $\text{hull}_\nabla(\text{gr}A)$ satisfy: $t \cdot \alpha(gtU^-) = \alpha(gU^-)$. Let $\Phi$: $\text{gr}A \to \text{hull}_\nabla(\text{gr}A)$ be defined by $\Phi(b)(gU^-) = \varphi(g^{-1} \cdot b)$ for all $g \in G$ and $b \in \text{gr}A$.

**Theorem 11** [G2; Theorems 8 and 16, p.131 and 133]. The mapping $\Phi$ is an injective, $G$-equivariant algebra homomorphism. It is surjective if and only if $A$ has a good filtration. The algebra $\text{hull}_\nabla(\text{gr}A)$ is integral over $\Phi(\text{gr}A)$ and $\Phi$ gives a surjective isomorphism of $(\text{gr}A)^U$ to $(\text{hull}_\nabla(\text{gr}A))^U$.

**Lemma 12.** Suppose that $A^U = k[a_1, \ldots, a_r]$ where $a_i$ has $T$-weight $\omega_i$. Suppose also that $\Phi(\text{gr}A)$ contains $a_i \otimes Y(\omega_i)$ for each $i = 1, \ldots, r$. Then, $\Phi$ is surjective.

**Proof.** Let $\omega, \omega' \in X^+(T)$. It is known that $Y(\omega)Y(\omega') = Y(\omega + \omega')$ [J; Proposition, p. 413]. Hence, for $a, a' \in A^U$ with $T$-weights $\omega, \omega'$, we have $(a \otimes Y(\omega))(a' \otimes Y(\omega')) = aa' \otimes Y(\omega + \omega')$. Since any $a \in A^U$ is a polynomial in the $a_i$'s, the lemma follows. □

**Note.** Later, we shall need the fact that $Y(\omega)Y(\omega') = Y(\omega + \omega')$ only when $G = \text{GL}_d$. Here, the proof is easily accessible by Young diagrams and straightening [G1; Theorem 13.5, p.79].

**Lemma 13.** Let $R$ and $S$ be a commutative $k$-algebras on which $G$ acts rationally. Suppose that $R \subset S$, $R^U = S^U$ and that $R$ has a good filtration. Then, $R = S$.

**Proof.** The mapping from $\text{gr}R$ to $\text{gr}S$ given by $(a + R_{n-1}) \to (a + S_{n-1})$ is an injective homomorphism so we may identify $\text{gr}R$ with a subalgebra of $\text{gr}S$. We have seen that there are homomorphisms $\Phi_R$: $\text{gr}R \to \text{hull}_\nabla(\text{gr}R) = (R^U \otimes k[G/U^-])^T$ and $\Phi_S$: $\text{gr}S \to \text{hull}_\nabla(\text{gr}S) = (S^U \otimes k[G/U^-])^T$. It follows from their definitions that the restriction of $\Phi_S$ to $\text{gr}R$ is $\Phi_R$. Also, by assumption $R^U = S^U$ so $\text{hull}_\nabla(\text{gr}R) = \text{hull}_\nabla(\text{gr}S)$.

    We now prove that $R_m = S_m$ by induction on $m$. First, $R_{-1} = S_{-1} = \{0\}$ by definition. In general, let $s \in S$ be in $S_{m+1}$, $s \notin S_m$. The mapping $\Phi_R$ is surjective by Theorem 11 since $R$ has a good filtration by assumption. Hence, there is an $r \in R_{i+1}$ with $\Phi_R(r + R_i) = \Phi_S(s + S_m)$. Then $\Phi_S(r + S_i) = \Phi_S(s + S_m)$. Since $\Phi_S$ is injective, $(r + S_i) = (s + S_m)$ so $i = m$ and $s - r$ is in $S_m = R_m$.

Thus, $s \in R$. □

We now apply Lemma 13 to the case where $G$ is $GL_d$, the subgroup $B$ consists of upper triangular matrices, the subgroup $U$ consists of those elements in $B$ with 1's on the diagonal, the subgroup $T$ consists of all diagonal matrices in $B$, $R = GL_d * j(^Hk[M_{n,n}])$ and $S = {}^Hk[M_{n,d}]$. By Lemma 2, $R^U = S^U$ so we have the following corollary to Lemma 13.

**Corollary 14.** If $GL_d * j(^Hk[M_{n,n}])$ has a good filtration, then $GL_d * j(^Hk[M_{n,n}]) = {}^Hk[M_{n,d}]$.

Since the actions of $H$ and $GL_d$ commute, $H$ sends the algebra $k[M_{n,d}]^U$ to itself. This gives an action of $H$ on $(k[M_{n,d}]^U \otimes k[G/U^-])^T$ when $H$ acts trivially on $k[G/U^-]$. Later, we shall need the following fact about this action.

**Lemma 15.** The mapping $\Phi$ is $H$-equivariant.

**Proof.** Since $H$ sends the algebra $k[M_{n,d}]^U$ to itself and leaves invariant each $T$-weight space, we have $\varphi(h \cdot b) = h \cdot \varphi(b)$ for all $b \in (gr k[M_{n,d}])^U$. Then,
$$\Phi(h \cdot b)(gU^-) = \varphi(g^{-1} * (h \cdot b)) = \varphi(h \cdot (g^{-1} * b)) = h \cdot \varphi(g^{-1} * b) = (h \cdot \Phi(b))(gU^-). \square$$

**Theorem 16.** Suppose that $^Hk[M_{n,d}]^U = k[a_1, \ldots, a_r]$ with $a_i$ having $T$-weight $\omega_i$ and that $\Phi(gr(GL_d * j(^Hk[M_{n,n}]))) \supset a_i \otimes Y(\omega_i)$. Then $GL_d * j(^Hk[M_{n,n}]) = {}^Hk[M_{n,d}]$.

**Proof.** We have seen (Lemma 2) that $^Hk[M_{n,d}]^U = (GL_d * j(^Hk[M_{n,n}]))^U$. Applying Lemma 12, we then see that $\Phi$ is surjective. According to Theorem 11, $GL_d * j(^Hk[M_{n,n}])$ has a good filtration and then is equal to $^Hk[M_{n,d}]$ by Corollary 14. □

**Corollary 17.** If $^Hk[M_{n,d}]$ is a finitely generated $k$-algebra and if $p$ is big enough, then $GL_d * j(^Hk[M_{n,n}]) = {}^Hk[M_{n,d}]$.

**Proof.** If $^Hk[M_{n,d}]$ is a finitely generated $k$-algebra, then so is $^Hk[M_{n,d}]^U$ [G1; Theorem 16.2, p.94]. Suppose that $^Hk[M_{n,d}]^U = k[a_1, \ldots, a_r]$ where $a_i$ has $T$-weight $\omega_i$. In general, if $a \in {}^Hk[M_{n,d}]^U$ has $T$-weight $\omega$ with $h(\omega) = n$, then $\Phi(a + {}^Hk[M_{n,d}]_{n-1}) = a \otimes f$ where $f$ is the unique element in $Y(\omega)$ which is fixed by $U$ and satisfies $f(e) = 1$ [G1; p.90]. If $p$ is sufficiently large, the representation $Y(\omega_i)$ is irreducible [J; Corollary 5.6, p.248]. Since $\Phi$ is injective and $G$ - equivariant, $\Phi(\mathrm{gr}(\mathrm{GL}_d * j(^H k[M_{n,n}]))) \supset a_i \otimes Y(\omega_i)$ and the hypotheses in Theorem 16 hold. □

## §5. Saturated subgroups

In this section, we obtain a partial converse to Corollary 14. We shall adapt the following notation for this section only. Let $A = k[M_{n,d}]$. Let $B_d$ be the subgroup of $\mathrm{GL}_d$ consisting of all upper triangular matrices, $U_d$ be the subgroup consisting of all those elements in $B_d$ with 1's on the diagonal and $T_d$ be the subgroup consisting of all diagonal matrices. Let $h_2: X(T_d) \to \mathbf{Z}$ be the homomorphism defined at the beginning of §4. Similarly, let $B_n$ be the subgroup of $\mathrm{GL}_n$ consisting of all upper triangular matrices, $U_n$ be the subgroup consisting of all those elements in $B_n$ with 1's on the diagonal and $T_n$ be the subgroup consisting of all diagonal matrices. Let $h_1: X(T_n) \to \mathbf{Z}$ be the homomorphism defined at the beginning of §4. Let $G = \mathrm{GL}_n \times \mathrm{GL}_d$, $B = B_n \times B_d$, $U = U_n \times U_d$ and $T = T_n \times T_d$. Let $A_m = \{a \in A : h_2(\chi_2) \le m$ for all weights $\chi = (\chi_1, \chi_2)$ of $T_d$ on $\langle \mathrm{GL}_d * a \rangle\}$. Since the actions of $\mathrm{GL}_n$ and $\mathrm{GL}_d$ commute, the group $\mathrm{GL}_n$ acts on each $A_m$.

**Theorem 18.** Each subspace $A_m$ has good $\mathrm{GL}_n$ - and $\mathrm{GL}_d$ - filtrations.

**Proof.** Let $m$ be a non-negative integer. The set $\pi_m = \{(\omega_1, \omega_2) \in X^+(T) : h_2(\omega_2) \le m\}$ is a saturated subset of $X^+(T)$. For let $(\mu_1, \mu_2), (\lambda_1, \lambda_2) \in X^+(T)$ with $(\mu_1, \mu_2) \le (\lambda_1, \lambda_2)$ and $(\lambda_1, \lambda_2) \in \pi_m$. Then $\lambda_2 - \mu_2 = \sum m_i \alpha_i$ where the $m_i$ are non-negative integers and the $\alpha_i$ are positive roots with respect to $B_d$. By property (ii),§4, of the homomorphism $h_2$, we have $h_2(\mu_2) \le h_2(\lambda_2) \le m$.

Let $O_{\pi, m}(A)$ be the maximal $G$ - submodule of $A$ belonging to $\pi_m$. We show now that $O_{\pi, m}(A) = A_m$. Indeed, let $a \in O_{\pi, m}(A)$. Then $\langle \mathrm{GL}_d * a \rangle \subset \langle G \cdot a \rangle \subset O_{\pi, m}(A)$ and, so, any highest weight $(\omega_1, \omega_2)$ of $G$ on $\langle G \cdot a \rangle$ satisfies $h_2(\omega_2) \le m$. Let $\chi_2$ be a weight of $T_d$ on $\langle \mathrm{GL}_d * a \rangle$. Since the actions of $\mathrm{GL}_n$ and $\mathrm{GL}_d$ commute, there is a vector $v \in \langle G \cdot a \rangle$ with $T$ - weight $\chi = (\chi_1, \chi_2)$ for some $\chi_1 \in X(T_n)$. Each weight $\omega$ of $T$ on $\langle U \cdot v \rangle$ satisfies $\chi \le \omega$. Choosing $\omega = (\omega_1, \omega_2)$ to be a highest weight, we see that $h_2(\chi_2) \le h_2(\omega_2) \le m$. Hence, $a \in A_m$. Next, let $a \in A_m$. The vector space $\langle G \cdot a \rangle = \langle \mathrm{GL}_n \cdot \langle \mathrm{GL}_d * a \rangle \rangle$. Also, if $v \in \langle \mathrm{GL}_d * a \rangle$ has $T_d$ weight $\chi_2$, then so does every element in $<GL_n \cdot v>$ since the actions of $GL_n$ and $GL_d$ commute. Since $a \in A_m$, $h_2(\chi_2) \leq m$. Hence, $h_2(\chi_2) \leq m$ for all weights $\chi = (\chi_1, \chi_2)$ of $T$ on $<G \cdot a>$ and, in particular, for any highest weight. Thus $a \in O_{\pi, m}(A)$ by the maximality of $O_{\pi, m}(A)$.

Now, $A$ has a good ($GL_n \times GL_d$) - filtration [D2; Proposition 1.3d, p.722]. Since $\pi_m$ is saturated, $O_{\pi, m}(A) = A_m$ has a good ($GL_n \times GL_d$) - filtration by (D1), §4. Then $A_m$ has a good $GL_n$ - filtration and $GL_d$ - filtration by [D2; Proposition 1.2e(i), p.720]. □

For $H$ a subgroup of $GL_n$, the following properties are equivalent: (i) $k[GL_n/H]$ has a good $GL_n$ - filtration where $GL_n$ acts by left translation; (ii) if $0 \to V_1 \to V_2 \to V_3 \to 0$ is a short exact sequence of $GL_n$ - modules where each $V_i$ has a good $GL_n$ - filtration, then the sequence $0 \to {}^H V_1 \to {}^H V_2 \to {}^H V_3 \to 0$ is exact. If $H$ satisfies either condition, we say that $H$ is a *saturated subgroup* of $GL_n$ [D1; p.121].

**Examples.** Any subgroup of $GL_n$ whose representations are completely reducible is a saturated subgroup by condition (ii). Thus, any finite group whose order is not divisible by $p$ or any torus is a saturated subgroup. We shall give other examples in §6, including the classical groups.

**Theorem 19.** If $H$ is a saturated subgroup of $GL_n$, then ${}^H k[M_{n,d}]$ has a good $GL_d$ - filtration.
**Proof.** The sequence $0 \to A_{m-1} \to A_m \to A_m/A_{m-1} \to 0$ is an exact sequence of $GL_n$ - modules each having a good $GL_n$ - filtration by Theorem 18 and (D2), §4. Since $H$ is a saturated subgroup of $GL_n$, the sequence $0 \to {}^H A_{m-1} \to {}^H A_m \to {}^H(A_m/A_{m-1}) \to 0$ is exact, i.e., ${}^H(A_m/A_{m-1}) = {}^H A_m/{}^H A_{m-1}$. Hence, ${}^H(grA) = gr({}^H A)$. Now, let $\Phi : grA \to \text{hull}_\nabla(grA)$ be as in §4. The mapping $\Phi$ is an isomorphism (by Theorem 11) since $A$ has a good $GL_d$ - filtration by Theorem 18 and (D2), §4. To show that ${}^H A$ has a good $GL_d$ filtration, we need to show that $\Phi : gr{}^H A \to \text{hull}_\nabla(gr{}^H A)$ is surjective (Theorem 11). By Lemma 15, $\Phi$ is $H$ - equivariant so ${}^H(\Phi(grA)) = \Phi({}^H(grA)) = \Phi(gr({}^H A))$. On the other hand, ${}^H(\Phi(grA))$ is the algebra of $H$ - invariants in $\text{hull}_\nabla(grA)$ since $\Phi(grA)) = \text{hull}_\nabla(grA)$. As defined just before Lemma 15, the action of $H$ on $\text{hull}_\nabla(grA) = (k[M_{n,d}]^U \otimes k[G/U^-])^T$ is the natural action on $k[M_{n,d}]^U$ and the trivial action on $k[G/U^-]$. Therefore, ${}^H(\text{hull}_\nabla(grA)) = ({}^H k[M_{n,d}]^U \otimes k[G/U^-])^T = \text{hull}_\nabla(gr{}^H A)$. □

**Corollary 20.** Let $H$ be a saturated subgroup of $GL_d$. Then $GL_d * j(^Hk[M_{n,n}]) = {}^Hk[M_{n,d}]$ if and only if $GL_d * j(^Hk[M_{n,n}])$ has a good $GL_d$-filtration.

**Proof.** Apply Corollary 14 and Theorem 19. □

## §6. Examples

### A. Finite groups

Before stating our main theorem here, we recall a theorem of Richman on vector invariants and then prove a lemma on linear functionals.

**Theorem 21** [R; Proposition 9, p.43]. Let $H$ be a finite subgroup of $GL_n$. Suppose that char $k = p$ and that $p$ divides $|H|$. Every set of $k$-algebra generators for $^Hk[M_{n,d}]$ contains a generator of degree $\geq d(p-1)/(p^{|H|-1}-1)$.

**Lemma 22.** Let $\ell$ be a linear functional on $M_{n,d}$. There is a linear functional $\ell'$ on $M_{n,n}$ and a $g \in GL_d$ so that $\ell = g * \ell'$.

**Proof.** Let $\ell$ be any linear functional on $M_{n,d}$, say $\ell = \sum c_{ij} x_{ij}$. We identify $\ell$ with the $n \times d$ matrix $(c_{ij})$. Then, a direct calculation shows that $g * \ell$ is identified with the matrix product $(c_{ij}) g^T$. Let $\text{Id}_n$ be the linear functional on $M_{n,n}$ corresponding to the $n \times d$ matrix which has the $n \times n$ identity matrix as its first $n$ columns and 0's elsewhere.

Now, fix a linear functional $\ell$ on $M_{n,d}$ and suppose that the $i$th row of $\ell$ is $v_i \in k^d$. If $\{v_1, \ldots, v_n\}$ is linearly independent, we choose $g$ to be any non-singular $d \times d$ matrix whose first $n$ rows are $v_1, \ldots, v_n$. Then, $\ell = g^T * \text{Id}_n$. Now, suppose that only $r$ rows of $\ell$ are linearly independent with $r < n$. We may assume that these rows are $v_1, \ldots, v_r$. Let $g$ be any non-singular $d \times d$ matrix whose first $r$ rows are $v_1, \ldots, v_r$. We construct a linear functional $\ell'$ on $M_{n,n}$ as follows. For $i = 1, \ldots, r$, row $i$ of $\ell'$ has 1 in the $(i, i)$ position and 0's elsewhere. For $i > r$, suppose that $v_i = a_1 v_1 + \ldots + a_r v_r$. Then, we take the $i$th row of $\ell'$ to be $a_1, \ldots, a_r, 0 \ldots, 0$. A calculation shows that $\ell = g^T * \ell'$. □

**Theorem 23.** Let $H$ be a finite subgroup of $\mathrm{GL}_n$ and let $p = \mathrm{char}\ k$.

(a) If $p > |H|$, then $\mathrm{GL}_d * j(^H k[\mathrm{M}_{n,n}]) = {}^H k[\mathrm{M}_{n,d}]$ and $^H k[\mathrm{M}_{n,d}]$ has a good filtration.

(b) If $p \mid |H|$, then $^H k[\mathrm{M}_{n,d}]$ is not equal to $\mathrm{GL}_d * j(^H k[\mathrm{M}_{n,n}])$ for all $d$ sufficiently large but is always its $p$-root closure.

(c) If $p < |H|$ and does not divide $|H|$, then $^H k[\mathrm{M}_{n,d}]$ may or may not equal $\mathrm{GL}_d * j(^H k[\mathrm{M}_{n,n}])$, but is always its $p$-root closure and always has a good $\mathrm{GL}_d$-filtration.

**Proof.** In statements (a) and (c), $p$ does not divide $|H|$ so the representations of $H$ are completely reducible. Then, condition (ii), §5, in the definition of saturated subgroup holds and we may apply Theorem 19 to see that $^H k[\mathrm{M}_{n,d}]$ has a good $\mathrm{GL}_d$-filtration. The assertions about $p$-root closure in statements (b) and (c) follow from Corollary 9. Now, we prove statement (a). Since $p > |H|$, $^H k[\mathrm{M}_{n,d}]$ is generated by orbit Chern classes, i.e., the various coefficients of $x$ in the polynomial $f = (x + h_1 \cdot \ell) \ldots (x + h_r \cdot \ell)$ where $\ell$ is a linear functional on $\mathrm{M}_{n,d}$ and $h_1 \cdot \ell, \ldots, h_r \cdot \ell$ are all the distinct elements in the $H$-orbit of $\ell$ [N - S; Theorem 4.1.2, p.80]. According to Lemma 22, there is a linear functional $\ell'$ on $\mathrm{M}_{n,n}$ and a $g \in \mathrm{GL}_d$ so that $\ell = g * \ell'$. Since the actions of $\mathrm{GL}_d$ and $H$ commute, the coefficients of $f = (x + h_1 \cdot (g * \ell')) \ldots (x + h_r \cdot (g * \ell'))$ are obtained by applying $g$ to the coefficients of $(x + h_1 \cdot \ell') \ldots (x + h_r \cdot \ell')$. But these coefficients are in $^H k[\mathrm{M}_{n,n}]$ so the coefficients of $f$ are in $\mathrm{GL}_d * j(^H k[\mathrm{M}_{n,n}])$.

To prove (b), we first note that if $\mathrm{GL}_d * j(^H k[\mathrm{M}_{n,n}]) = {}^H k[\mathrm{M}_{n,d}]$, then there is a positive integer $N$ so that $^H k[\mathrm{M}_{n,d}]$ is generated by polynomials of degree $\leq N$ for all $d$. In fact, suppose that $^H k[\mathrm{M}_{n,n}]$ is generated by $f_1, \ldots, f_r$. Then, $\mathrm{GL}_d * j(^H k[\mathrm{M}_{n,n}])$ is generated by bases for the vector spaces $<\mathrm{GL}_d * f_1>, \ldots, <\mathrm{GL}_d * f_r>$ and we may take $N$ to be the maximum degree of the $f_i$. Now, Theorem 21 says that there is no degree bound on generators of the $^H k[\mathrm{M}_{n,d}]$ for all $d$. Hence, $\mathrm{GL}_d * j(^H k[\mathrm{M}_{n,n}])$ cannot be $^H k[\mathrm{M}_{n,d}]$ for sufficiently large $d$.

As for statement (c), examples show that $^H k[\mathrm{M}_{n,d}]$ may or may not equal $\mathrm{GL}_d * j(^H k[\mathrm{M}_{n,n}])$. For instance, let $\mathrm{char}\ k = 2$ and let $H$ be the group consisting of the three cube roots of unity in $k$, say $\{a, a^2, 1\}$. Let $\rho$ be the representation of $H$ which sends $a$ to the $2 \times 2$ diagonal matrix with $a, a$ on the diagonal. We shall denote entries in the first row of $m \in \mathrm{M}_{2,d}$ by $x_1, \ldots, x_d$ and in the second row by $y_1, \ldots, y_d$. Then, $^H k[\mathrm{M}_{2,d}]$ is generated by all monomials of degree 3 and is not $\mathrm{GL}_d * j(^H k[\mathrm{M}_{2,2}])$ since $x_1 x_2 x_3$ is not in $\mathrm{GL}_d * j(^H k[\mathrm{M}_{2,2}])$. Next, let $\rho$ be the representation of $H$ which sends $a$ to the $3 \times 3$ diagonal matrix with $a$, $a$, $a$ on the diagonal. Then, $^H k[M_{3,d}]$ is generated by all monomials of degree 3 and is equal to $GL_d * j(^H k[M_{3,3}])$. □

**Notes.** With respect to (c), Knop showed that $^H k[M_{n,d}]$ can be obtained by polarizing the invariants of $d_o$ vectors where $d_o \geq \max\{n, \beta_k(H)/(p-1)\}$ [K]. Next, the assertions about $p$ - root closure in Theorem 23 can be proved in a different way. Let $X$ (resp. $Y$) be the affine variety such that $k[X] = {}^H k[M_{n,d}]$ (resp. $k[Y] = GL_d * j(^H k[M_{n,n}])$). Let $\pi: X \to Y$ be the map corresponding to the inclusion of $k[Y]$ in $k[X]$; $\pi$ is proper and surjective by Theorem 3. It is injective since $GL_d * j(^H k[M_{n,n}]))$ separates $H$ - orbits by the results in [DKW]. Thus, we may apply Theorem 6 to get the assertions on $p$ -root closure.

**B. Classical groups**

Our purpose here is to derive the well - known facts about vector invariants of the classical groups in arbitrary characteristic from the results proved in §4. Before doing that, however, we prove a general lemma and then give a procedure for simplifying our calculations. With these tools in place, the facts about the classical groups will follow quickly.

**Lemma 24.** Let $H$ be a subgroup of $GL_n$.
(a) Let $R$ be a subalgebra of $GL_d * j(^H k[M_{n,n}])$ on which $GL_d$ acts and which satisfies the following two conditions: (1) if $^H k[M_{n,d}]^U = k[a_1, \ldots, a_r]$ where $a_i$ has $T$ - weight $\omega_i$, then each $a_i$ is in $R$; (2) $\Phi(\text{gr} R)$ contains $a_i \otimes Y(\omega_i)$ for $i = 1, \ldots, r$. Then $R = GL_d * j(^H k[M_{n,n}]) = {}^H k[M_{n,d}]$.
(b) Suppose that $H$ is a subgroup of $SL_n$ and that $^H k[M_{n,n}]^U = k[a_1, \ldots, a_r]$ where $a_i$ has $T$ - weight $\omega_i$. If $\Phi(\text{gr}(^H k[M_{n,n}]))$ contains $a_i \otimes Y(\omega_i)$ for $i = 1, \ldots, r$, then $H$ is a saturated subgroup of $GL_n$.
**Proof.** We first prove statement (a). According to (1), $R^U = {}^H k[M_{n,d}]^U = k[a_1, \ldots, a_r]$. Since (2) holds, we may apply Lemma 12 to see that $\Phi: \text{gr} R \to \text{hull}_\nabla(\text{gr} R)$ is surjective. Then, $R$ has a good filtration by Theorem 11 and, applying Lemma 13, we see that $R = {}^H k[M_{n,d}]$ since $R^U = {}^H k[M_{n,d}]^U$. Since $R$ is contained in $GL_d * j(^H k[M_{n,n}])$, we see that $GL_d * j(^H k[M_{n,n}]) = {}^H k[M_{n,d}]$.

To prove statement (b), we need to show that the action of $GL_n$ by right translation on

$^H k[GL_n]$ admits a good filtration. Since $H$ is a subgroup of $SL_n$, $^H k[GL_n] = {}^H k[M_{n,n}][1/\Delta]$ and so $^H k[GL_n]^U = {}^H k[M_{n,n}]^U [1/\Delta] = k[a_1, \ldots, a_r][1/\Delta]$. Now, $\Phi(\text{gr}(^H k[GL_n])) \supset \Phi(\text{gr}(^H k[M_{n,n}]))$ and, by assumption, $\Phi(\text{gr}(^H k[M_{n,n}]))$ contains $a_i \otimes Y(\omega_i)$ for each $i = 1, \ldots, r$  A calculation shows that $\Phi(1/\Delta) = (1/\Delta) \otimes \det(g)$. Hence, $\Phi(\text{gr}(^H k[GL_n])) \to \text{hull}_\nabla(\text{gr}(^H k[GL_n]))$ is surjective by Lemma 12 and $H$ is a saturated subgroup of $GL_n$ by Theorem 11. □

The calculation in (2) is simplified if we replace the mapping $\Phi$ by another mapping $\Phi'$ which we now describe. For a moment, let $G$ be any reductive group and let $A$ be a commutative $k$-algebra on which $G$ acts rationally via a mapping $G \times A \to A$ denoted by $(g, a) \to g \cdot a$. We shall use the notation introduced at the beginning of §4. The algebra $k[G/U^-]$ is isomorphic to the algebra $k[U^-\backslash G]$ via the mapping $\epsilon: k[G] \to k[G]$ given by $(\epsilon f)(g) = f(g^{-1})$. The group $T$ acts on $k[G/U^-]$ by right translation and on $k[U^-\backslash G]$ by left translation and $\epsilon$ is $T$-equivariant. The mapping $\text{id} \otimes \epsilon$ then gives an isomorphism from $\text{hull}_\nabla(\text{gr}A) = (A^U \otimes_k k[G/U^-])^T$ to $(A^U \otimes_k k[U^-\backslash G])^T$. The injective homomorphism $\Phi: \text{gr}A \to (A^U \otimes k[G/U^-])^T$ gives rise to an injective homomorphism $\Phi': \text{gr}A \to (A^U \otimes k[U^-\backslash G])^T$ by $\Phi' = (\text{id} \otimes \epsilon)\Phi$. Thus, $\Phi'(b)(U^- g) = \Phi(b)(g^{-1}U^-) = \varphi(g \cdot b)$. Obviously, $\Phi$ is surjective if and only if $\Phi'$ is surjective.

For $\omega \in X^+(T)$, we define subspaces $Y'(\omega)$ of $k[U^-\backslash G]$ as follows: $Y'(\omega) = \{f \in k[U^-\backslash G] : \ell_t f = \omega(t)^{-1} f \text{ for all } t \in T\}$. The group $G$ acts on $Y'(\omega)$ by right translation and $\epsilon$ is a $G$-equivariant isomorphism from $Y(\omega)$ to $Y'(\omega)$. If we choose a basis $\{a_i\}$ of $A^U$ so that $t \cdot a_i = \omega_i(t) a_i$ for all $t \in T$, then we see that $(A^U \otimes k[U^-\backslash G])^T$ is a direct sum of the subspaces $a_i \otimes Y'(\omega_i)$.

When $G = GL_d$ and $B = TU$ is the usual subgroup of upper triangular matrices, the subspaces $Y'(\omega)$ can be described as follows [G1; Theorem 13.5, p.79]. Let $\omega = e_1\omega_1 + \ldots + e_d\omega_d$ where $e_d \geq 0$. A basis of $Y'(\omega)$ consists of all standard bitableaux $(D \mid E)$ where each row of $D$ has the form $1\ 2\ \ldots\ r$ and there are $e_1$ rows of length 1, $e_2$ rows of length 2, $\ldots$, $e_d$ rows of length $d$.

Let $(i_1 \ldots i_r \mid j_1 \ldots j_r)_m$ denote the minor of $m \in M_{n,d}$ corresponding to the choice of $\{i_1 \ldots i_r\}$ as rows and $\{j_1 \ldots j_r\}$ as columns. We write $(i_1 \ldots i_r \mid j_1 \ldots j_r)_g$ for the corresponding minor of $g \in GL_d$.

**Theorem 25.** (a) $SL_n$ is a saturated subgroup of $GL_n$; furthermore, the algebra of $SL_n$ - invariants on $M_{n,d}$ is generated by all $n \times n$ minors.

(b) If char $k \neq 2$, then $SO_n$ is a saturated subgroup of $GL_n$. Furthermore, the algebra of $SO_n$ - invariants on $M_{n,d}$ is generated by all $n \times n$ minors and all inner products $\langle x_i, x_j \rangle$.

(c) If char $k \neq 2$, then $Sp_{2n}$ is a saturated subgroup of $GL_{2n}$. Furthermore, the algebra of $Sp_{2n}$ - invariants on $M_{2n,d}$ is generated by all skew-symmetric products $\langle x_i, x_j \rangle$.

**Proof.** We sketch the proof for (b). This proof also contains the proof of (a). The proof of (c) is similar when the determinant is replaced by the Pfaffian.

We show that the conditions in Lemma 24 hold for all $d \geq n$ when $H = SO_n$ and char $k \neq 2$. Let $R$ be the subalgebra of $^H k[M_{n,d}]$ generated by all $n \times n$ minors and all inner products $\langle x_i, x_j \rangle = \sum x_{qi} x_{qj}$. The inner product $\langle x_i, x_j \rangle$ is in $GL_d * j(^H k[M_{n,n}])$ since it can be written as $g * \langle x_1, x_2 \rangle$ where $g$ is a suitable permutation matrix in $GL_d$. The $n \times n$ minor $(1\,2 \ldots n \mid i_1 \ldots i_n)_m$ is also in $GL_d * j(^H k[M_{n,n}])$. Indeed,

$$g * \Delta = \sum (1\,2 \ldots n \mid i_1 \ldots i_n)_m (i_1 \ldots i_n \mid 1\,2 \ldots n)_g.$$

Each term appearing in this sum has a different $T$-weight and so is in $\langle GL_d * \Delta \rangle$. That $GL_d$ sends $R$ to itself is shown by the following two formulas.

(i) $g * \langle x_i, x_j \rangle = \sum \langle x_q, x_t \rangle g_{qi} g_{tj}$.

(ii) $g * (1\,2 \ldots n \mid j_1 \ldots j_n)_m = \sum (1\,2 \ldots n \mid i_1 \ldots i_n)_m (i_1 \ldots i_n \mid j_1 \ldots j_n)_g$.

Next, we check statements (1) and (2) in Lemma 24(a). The algebra of $U$ invariants on $^H k[M_{n,d}]$ is $k[D_1, \ldots, D_{n-1}, \Delta]$ where $D_r$ is the determinant of the $r \times r$ matrix whose $(i,j)$ entry is $\langle x_i, x_j \rangle$ for $i, j \in \{1, 2, \ldots, r\}$; $D_r$ has $T$-weight $2\omega_r$ and is in $R$ [G1; Theorem 11.9, p. 64]. Let $\langle i_1 \ldots i_r \mid j_1 \ldots j_r \rangle$ denote the determinant of the $r \times r$ matrix whose $(h,q)$-entry is $\langle x_h, x_q \rangle$ for $h \in \{i_1 \ldots i_r\}$ and $q \in \{j_1 \ldots j_r\}$. A direct calculation shows that

(iii) $\Phi'((1\,2 \ldots n \mid j_1 \ldots j_n)_m + R_{n(d-n)-1}) = \Delta \otimes (1\,2 \ldots n \mid j_1 \ldots j_n)_g$ and

(iv) $\Phi'(\langle i_1 \ldots i_r \mid j_1 \ldots j_r \rangle + R_{2r(d-r)-1}) = D_r \otimes [(1\,2 \ldots r \mid i_1 \ldots i_r)_g (1\,2 \ldots r \mid j_1 \ldots j_r)_g]$.

That $\Phi'(\text{gr}R)$ contains $a_i \otimes Y'(\omega_i)$ for each $i = 1, \ldots, r$ now follows from (iii) and (iv) and the explicit description of $Y'(\omega)$ given above. $\square$

## C. Certain regular unipotent subgroups

We begin with a general theorem and then apply it to a class of regular unipotent subgroups discovered by Pommerening.

**Theorem 26.** Let $G$ be an algebraic subgroup of $GL_n$ and let $H$ be a closed subgroup of $G$ such that $k[G/H]$ is a finitely generated $k$-algebra. Let $Z$ be the affine variety such that $k[Z] = k[G/H]$. Suppose that there is a point $z_o$ in $Z$ so that $gz_o = z_o$ for all $g \in G$. If $GL_d * j(^H k[M_{n,n}]) = {}^H k[M_{n,d}]$, then $GL_d * j(^G k[M_{n,n}]) = {}^G k[M_{n,d}]$.

**Proof.** For every $d \geq n$, we define two $k$-algebra mappings, $\rho_d$ and $\psi_d$. The mapping $\rho_d : {}^G(k[M_{n,d}] \otimes k[G/H]) \to {}^H k[M_{n,d}]$ given by $\rho_d(\sum a_i \otimes f_i) = \sum f_i(e) a_i$ is a surjective isomorphism [G1, Theorems 4.3, p.20 and 9.1, p. 49 or Exercise 3, p. 53]. Here, $G$ acts on $k[M_{n,d}]$ and $k[G/H]$ by left translation. Next, the mapping $\psi_d : {}^G(k[M_{n,d}] \otimes k[G/H]) \to {}^G k[M_{n,d}]$ given by $\psi_d(\sum a_i \otimes f_i) = \sum f_i(z_o) a_i$ is a surjective homomorphism [BK; Corollary 3.7, p.89].

The group $GL_d$ acts on $k[M_{n,d}] \otimes k[G/H]$ by $g * (\sum a_i \otimes f_i) = \sum (g * a_i) \otimes f_i$. Since the actions of $G$ and $GL_d$ commute, $GL_d$ sends ${}^G(k[M_{n,d}] \otimes k[G/H])$ to itself. A short computation shows that both $\rho_d$ and $\psi_d$ are $GL_d$-equivariant. Furthermore, we have $j \circ \psi_n \circ \rho_n^{-1} = \psi_d \circ \rho_d^{-1} \circ j$. This follows from the fact that if $f \in {}^H k[M_{n,n}]$ is $\rho_n(\sum a_i \otimes f_i)$ where $a_i \in k[M_{n,n}]$, then $jf = \rho_d(\sum ja_i \otimes f_i)$. Then,

$${}^G k[M_{n,d}] = \psi_d({}^G(k[M_{n,d}] \otimes k[G/H])) = (\psi_d \circ \rho_d^{-1})({}^H k[M_{n,d}]) = (\psi_d \circ \rho_d^{-1})(GL_d * j(^H k[M_{n,n}]))$$
$$= GL_d * j((\psi_n \circ \rho_n^{-1})(^H k[M_{n,n}])) = GL_d * j(^G k[M_{n,n}])). \quad \square$$

In the rest of this section, we use facts about observable groups. Let $T_n$ be the maximal torus in $GL_n$ consisting of all diagonal matrices with non-zero determinant. Let $H$ be a subgroup of $GL_n$ which is normalized by $T_n$. We shall denote the observable hull of $T_n H$ by $(T_n H)''$; it is the smallest reductive group in $GL_n$ which contains $T_n H$ [G1; Lemma 3.10, p.18]. Furthermore, if $V$ is a rational $GL_n$-module and if $T_n H$ fixes a point $v$ in $V$, then so does $(T_n H)''$ [G1; p.6].

**Corollary 27.** Let $H$ be a unipotent subgroup of $GL_n$ which is normalized by $T_n$ and is such that $k[GL_n/H]$ is a finitely generated $k$-algebra. Let $L$ be any algebraic subgroup of $GL_n$ such that $H \subset L \subset (T_n H)''$. If $GL_d * j(^H k[M_{n,n}]) = {}^H k[M_{n,d}]$, then $GL_d * j(^L k[M_{n,n}]) = {}^L k[M_{n,d}]$. If $H$ is a saturated subgroup of $GL_n$, then so is $L$.

**Proof.** First, $k[L/H]$ is a finitely generated $k$ - algebra since $k[GL_n/H]$ is [G1; Corollary 4.5, p.22]. Now, $k[L/H] = k \oplus M$ where $M$ is the sum of all non-zero $T_n$ weight spaces where $T_n$ acts by right translation. Since $L \subset (T_nH)''$, $M$ is an $L$ - invariant ideal. Therefore, the assumptions in Theorem 26 hold with $G$ replaced by $L$. Next, suppose that $H$ is a saturated subgroup of $GL_n$ and let $0 \to V_1 \to V_2 \to V_3 \to 0$ is a short exact sequence of $GL_n$ - modules where each $V_i$ has a good $GL_n$ - filtration. By assumption, the sequence $0 \to {}^HV_1 \to {}^HV_2 \to {}^HV_3 \to 0$ is exact. The representations of $T_n$ are completely reducible so the sequence $0 \to {}^{T_nH}V_1 \to {}^{T_nH}V_2 \to {}^{T_nH}V_3 \to 0$ is also exact. Since $L \subset (T_nH)''$, we see that the sequence $0 \to {}^LV_1 \to {}^LV_2 \to {}^LV_3 \to 0$ is exact. □

Pommerening described a class of unipotent subgroups $H$ in $GL_n$ for which the assumptions in Corollary 27 hold [P]. This class contains the unipotent radicals of parabolic subgroups of $GL_n$. Any subgroup $H$ in this class has the property that $R = {}^Hk[M_{n,d}]$ is spanned by standard bitableaux, each row of which is an invariant $H$ - minor [P; p.277]. Any such minor is obtained by applying a column permutation to an $H$ - invariant minor in $k[M_{n,n}]$. Thus, $GL_d * j({}^Hk[M_{n,n}]) = {}^Hk[M_{n,d}]$. Also, it follows that ${}^Hk[M_{n,d}]^U$ is generated by $H$ - invariant minors of the form $(i_1 \ldots i_r | 1\ 2 \ldots r)$. Since $\Phi'((i_1 \ldots i_r | j_1 \ldots j_r)_m + R_{r(d-r)-1}) = (i_1 \ldots i_r | 1\ 2 \ldots r)_m \otimes (1\ 2 \ldots r | j_1 \ldots j_r)_g$, we see that any such $H$ is a saturated subgroup of $GL_n$.